\documentclass[a4paper,12pt]{article}
\usepackage{amsthm}
\usepackage{amssymb}
\usepackage{amsmath}
\usepackage{amsfonts}

\newtheorem{lemma}{Lemma}[section]
\newtheorem{theorem}{Theorem}[section]

\newtheorem{definition}{Definition}[section]

\def\b1{\mbox{\boldmath $1$}}

\newenvironment{demo*}{\vspace{3mm}\noindent{\bf Proof.}}{\hfill $\Box$ \vspace{3mm}}

\begin{document}

\baselineskip=13pt

\title{\bf  On optimality of the barrier strategy for a general L\'evy risk process}

\author{\normalsize$^a$KAM CHUEN YUEN,  and $^b$CHUANCUN YIN\\
{\normalsize\it $^a$Department of Statistics and Actuarial Science,
The University of Hong Kong,} \\
\noindent{\normalsize\it Pokfulam Road, Hong Kong}\\ e-mail: kcyuen@hku.hk\\
[3mm]
{\normalsize\it $^b$School of Mathematical Sciences, Qufu Normal University}\\
{\normalsize\it Shandong 273165, P.R.\ China}\\
e-mail: ccyin@mail.qfnu.edu.cn }
\date{}
\maketitle

 \noindent{\large {\bf Abstract}} \ We consider the optimal dividend problem for the insurance risk process in a general L\'evy process setting. The objective is to find a strategy which maximizes the expected total discounted dividends until the time of ruin. We give sufficient conditions under which the optimal strategy is of barrier type. In particular, we show that if the L\'evy density is a completely monotone function, then the optimal dividend  strategy   is   a barrier strategy.  This approach was inspired by the work of  Avram et al. (2007) [Annals of Applied Probability 17, 156-180],  Loeffen (2008) [Annals of Applied Probability 18, 1669-1680] and Kyprianou et al. (2010) [Journal of Theoretical Probability 23, 547-564]
 in which the same problem was considered under the spectrally negative L\'evy processes
 setting.

\noindent {\bf Keywords}\;  { L\'evy processes $\cdot$ Optimal
dividend problem  $\cdot$ Complete monotonicity   $\cdot$ Barrier
strategy $\cdot$ Scale function  $\cdot$  Probability of ruin

}


\normalsize

\baselineskip=13pt

\section{The model and problem setting}\label{intro}
Let  $X = \{X_t : t \ge 0\}$ be  a real-valued  L\'evy process
define on a filtered probability space $(\Omega, \cal{F}, \Bbb{F},
P)$ where $\Bbb{F}=(\cal{F})$$_{t\ge 0}$ is generated by the process
$X$ and  satisfies the usual conditions of right-continuity and
completeness. Denote by $P_x$ the probability law of $X$ when it
starts at $x$. For notational convenience, we write $P=P_0$. Let
$E_x$ be the expectation operator associated with  $P_x$ with
$E=E_0$. For $\theta\in \Bbb{R}$, $\kappa(\theta)$ denotes the
characteristic exponent of $X$ given by
$$\kappa(\theta)=\frac{1}{t}\log E(e^{i\theta X_t}) =ia\theta-\frac1 2\sigma^2\theta^2
+\int_{-\infty}^{\infty}(e^{i\theta x}-1-i\theta x\text{\bf
1}_{\{|x|<1\}})\Pi (dx),$$ where $a,\sigma$ are real constants, and
$\Pi$ is a positive measure supported  on $(-\infty,\infty)\setminus
\{0\}$ which satisfies the integrability condition
$$\int_{-\infty}^{\infty}\min\{1,x^2\}\Pi(dx)<\infty.$$
The characteristics $(a, \sigma^2, \Pi)$ are called the L\'evy
triplet of the process and completely determines its law; $\Pi$ is
called the L\'evy measure; and $\sigma$ is the Gaussian component of
$X$. If $\Pi(dx) =\pi(x)dx$, then we  call $\pi$  the L\'evy
density. Such a L\'evy process is of  bounded variation if and only
if $\sigma=0$ and
$\int_{-\infty}^{\infty}\min\{|x|,1\}\Pi(dx)<\infty$.  It is well
known that L\'evy process $X$ is a space-homogeneous strong Markov
process.

For $\text{{\bf Re}}\theta=0$,  we can define the Laplace exponent
of the process X by
\begin{equation}
  \Psi (\theta)=\frac1{t}\log E(e^{\theta X_t}) =a\theta + \frac1 2\sigma^2\theta^2
+\int_{-\infty}^{\infty}(e^{\theta x}-1-\theta x\text{\bf
1}_{\{|x|<1\}})\Pi (dx). \label{intro-eq1}
\end{equation}
That is,
$$E(e^{\theta X_t}) = e^{t\Psi(\theta)}, \quad \text{{\bf Re}}\theta=0, \quad t\ge 0.$$
Using Ito's formula, we find that the infinitesimal generator for
$X$ is given by
\begin{equation}
\Gamma g(x)=\frac{1} {2}\sigma^2 g''(x)+a g'(x)
+\int_{-\infty}^{\infty}[g(x+y)-g(x)-g'(x)y\text{\bf
1}_{\{|y|<1\}}]\Pi(dy), \label{intro-eq2}
\end{equation}
for $g\in C^2$ with compact support, where $\text{\bf 1}_{A}$ is the
indicator function of set $A$.

When $\Pi\{(0,\infty)\}=0$, i.e., the L\'evy process $X$ with no
positive jumps, is called the spectrally negative L\'evy process. In
this case, we recall from  Bertoin (1998) and Kyprianou (2006) that
for each $q\ge 0$, there exists a continuous and increasing function
$W^{(q)}:\Bbb{R}\rightarrow [0,\infty)$, called the $q$-scale
function defined in such a way that $W^{(q)}(x) = 0$ for all $x < 0$
and on $[0,\infty)$ its Laplace transform is given by
$$\int_0^{\infty}\text{e}^{-\theta
x}W^{(q)}(x)dx=\frac{1}{\Psi(\theta)-q},\; \theta >\phi(q),$$ where
$\phi(q)=\sup\{\theta\ge 0:\Psi(\theta)=q\}$ is the right inverse of
$\Psi$. Smoothness of the scale function is related to the
smoothness of the underlying paths of the associated process. The
following facts are taken from Kyprianou et al. (2010) and Chan et
al. (2010). It is known that if X has paths of bounded variation
then, for all $q\ge 0$, $W^{(q)}|_{(0,\infty)}\in C^1(0,\infty)$ if
and only if $\Pi$ has no atoms. In the case that $X$ has paths of
unbounded variation, it is known that, for all $q\ge 0$,
$W^{(q)}|_{(0,\infty)}\in C^1(0,\infty)$. Moreover, if $\sigma> 0$,
then $C^1(0,\infty)$ may be replaced by $C^2(0,\infty)$; and if the
L\'evy measure has a density, then the scale function is always
differentiable. In particular, if $\pi$ is completely monotone, then
$W^{(q)}|_{(0,\infty)}\in C^{\infty}(0,\infty)$.

For general references on L\'evy processes, we refer the reader,
among others, to Bertoin (1998), Sato (1999) and Kyprianou (2006).

In this paper, we only consider the case that $\Pi$ with density
$\pi$ is absolutely continuous with respect to Lebesgue measure. For
the L\'evy process $X$, we consider the following de Finetti's
dividend problem. Let $ \xi=\{L_t^{ \xi}:t\ge 0\}$ be a dividend
strategy consisting of a left-continuous non-negative non-decreasing
process adapted to the filtration $\{\cal{F}$$_t$\}$_{t\ge 0}$ of
$X$. Specifically, $L_t^{ \xi}$ represents the cumulative dividends
paid out up to time $t$ under the control $\xi$ for an insurance
company whose risk process is modelled by $X$. We define the
controlled risk process $U^{\xi}=\{U_t^{\xi}:t\ge 0\}$ by $U_t^{
\xi}=X_t-L_t^{ \xi}$. Let $\tau^{\xi}=\inf\{t>0: U_t^{ \xi}<0\}$ be
the ruin time when the dividend payments are taken into account.
Define the value function of a dividend strategy $ \xi$ by
$$V_{\xi}(x)=E_x\left(\int_0^{\tau^{ \xi}}\text{e}^{-\delta t}dL_t^{
\xi}\right),$$ where $\delta>0$ is the discounted rate. The integral
is understood pathwise in a Lebesgue-Stieltjes sense.

A dividend strategy is called admissible if $L_{t+}^{ \xi}-L_t^{
\xi}\le U_t^{ \xi}$ for $t<\tau^{ \xi}$. In words, the lump sum
dividend payment is smaller than the size of the available capitals.
Let $\Xi$ be the set of all admissible dividend policies. The de
Finetti dividend problem consists of solving the following
stochastic control problem:
$$V_*(x)=\sup_{ \xi\in\Xi}V_{ \xi}(x),$$
and, if it exists, we want to find a strategy $ \xi_*\in\Xi$ such
that $V_{ \xi_*}(x)=V_*(x)$ for all $x\ge 0$.

This optimization problem goes back to de Finetti (1957), who
considered a discrete time random walk with step sizes $\pm 1$ and
showed that a certain barrier strategy maximizes expected discounted
dividend payments. Optimal dividend problem has recently gained
great attention in the actuarial literature.  For the diffusion risk
process, the optimal problem has been studied by many authors
including Asmussen et al. (2000), Paulsen (2003), and Decamps and
Villeneuve (2007). It is well known that under some reasonable
assumptions, the optimality in the diffusion process setting is
achieved by a barrier strategy (see, for example, Shreve et al.
(1984)). The general problem for the Cram\'er-Lundberg risk model
was first solved by Gerber in (1969) via a limit of an associated
discrete problem. Recently, Azcue and Muler (2005) used the
technique of stochastic control theory and Hamilton-Jacobi-Bellman
(HJB) equation to solve the problem. They also included a general
reinsurance strategy as another possible control.  For the
Cram\'er-Lundberg risk model with interest, Yuen et al. (2007)
investigated some ruin problems in the presence of a constant
dividend barrier; and Albrecher and Thonhauser (2008) derived the
optimal dividend strategy which is again of band type and for
exponential claim sizes collapses to a barrier strategy. In fact,
for various risk models, many results in the literature indicate
that a band strategy turns out to be optimal among all admissible
strategies.  For a more general risk process, namely the spectrally
negative L\'evy process, Avram et al. (2007) gave a sufficient
condition involving the generator of the L\'evy process for the
optimality of the barrier strategy; Loeffen (2008) connected the
shape of the scale function to the existence of an optimal barrier
strategy and showed that the optimal strategy is a barrier strategy
if the L\'evy measure has a completely monotone density; and
Kyprianou et al. (2010) further investigated the optimal dividend
control problem and showed that the problem is solved by a barrier
strategy whenever the L\'evy measure of a spectrally negative L\'evy
process has a density which is log-convex.

Motivated by the work of Avram et al. (2007), Loeffen (2008) and
Kyprianou et al. (2010) for a spectrally negative L\'evy process,
our objective is to consider the  optimal dividend problem for a
general L\'evy process (not necessarily spectrally negative). In
contrast to the previous works, our approach does not rely on the
theory of scale function  of a spectrally negative L\'evy process.
Instead, the derivation of our results requires an introduction of a
new function and the use of the Wiener-Hopf factorization theory.
The rest of this paper is organized as follows. In Section 2, we
present some preliminary results which are derived in Rogers (1983),
Bertion and Doney (1994) and Kyprianou (2006). In Section 3, we
discuss a barrier strategy for dividend payments for the risk model
of study. In Section 4, we give the main results and their proofs.
Finally, Section 5 presents some examples and ends with a remark.

\vskip 0.2cm
\section{Two lemmas on the probability of ruin }\label{con}
\setcounter{equation}{0}

In this section, we present two results on the probability of ruin
which will be used later.

The first result is about the complete monotonicity of the
probability of ruin. By definition, an infinitely differentiable
function $f\in (0,\infty)\rightarrow [0,\infty)$  is called
completely monotone if $(-1)^n f^{(n)}(x)\ge 0$ for all
$n=0,1,2,\cdots$.  Denote by $\tau(q)$ an exponential random
variable with mean $1/q$ which is independent of the process $X$.
For $q=0$, $\tau(0)$ is understood to be infinite. Furthermore, let
$$\underline{X}_{\tau(q)}=\inf_{0\le t\le\tau(q)}X_t \quad \text{and}  \quad \overline{X}_{\tau(q)}=\sup_{0\le t\le\tau(q)}X_t,$$
be the infimum and the supremum of the the L\'evy process $X$ killed
at the random time $\tau(q)$, respectively. In the case with $q=0$,
we always assume that
\begin{equation}
E(X_1)=\Psi'(0) =a +\int_{-\infty}^{\infty}x\text{\bf 1}_{\{|x|\ge
1\}}\pi(x)dx>0.\label{con-eq1}
\end{equation}

We call the characteristic functions $E(\exp(\alpha
\underline{X}_{\tau(q)}))$ and $E(\exp(\alpha
\overline{X}_{\tau(q)}))$ the right and left Wiener-Hopf factors of
$X$, respectively. According to Rogers (1983), the Wiener-Hopf
factors are called mixtures of exponentials if there are probability
measures $H_{+}$ and $H_-$ on $(0,\infty]$ such that
$$E(\exp(\alpha \overline{X}_{\tau(q)}))=\int_0^{\infty} \frac{\lambda H_{+}(d\lambda)}{\lambda-\alpha},
 \quad E(\exp(\alpha \underline{X}_{\tau(q)}))=\int_0^{\infty}
\frac{ \lambda H_{-}(d\lambda)}{\lambda+\alpha}.$$ Note that if
$E\exp(\alpha \underline{X}_{\tau(q)})$ is a mixture of
exponentials, then the probability distribution function of
$\underline{X}_{\tau(q)}$ has the form
$$P_{-}(x)=
\int_{(0,\infty]}e^{\lambda x}H_{-}(d\lambda), \quad x<0.$$
Similarly,
  the probability  distribution function  of $\overline{X}_{\tau(q)}$ has the form
$$P_{+}(x)=
\int_{(0,\infty]}(1-e^{-\lambda x})H_{+}(d\lambda), \quad x>0.$$

For $x\ge 0$, let $\tau=\inf\{t\ge 0: x+X(t)\le 0\}$. Define its
Laplace transform and the probability of ruin as
$$\psi_q(x)=E(e^{-q\tau}), \quad \psi(x)=P(\tau<\infty),$$
respectively. Then,
$$\psi_q(x)=P(\underline{X}_{\tau(q)}\le -x),$$
and
$$\psi(x)=\lim_{q\to 0}\psi_q(x)=P(\underline{X}_{\infty}\le
-x).$$

\begin{definition} (Rogers (1983))\; The L\'evy process $X$ has completely monotone
L\'evy density if there exist measures $\mu_{+}, \mu_{-}$ on
$(0,\infty)$ such that
\begin{equation}
\pi(x)=\text{\bf 1}_{\{x>0\}}\int_{(0,\infty)}e^{-tx}\mu_{+}(dt)
+\text{\bf 1}_{\{x<0\}}\int_{(0,\infty)}e^{tx}\mu_{-}(dt),
\label{con-eq2}
\end{equation}
 where
$$\int \frac{1}{t(1+t)^2}(\mu_{+}+\mu_{-})(dt)<\infty.$$
\end{definition}

Theorem 2 of Rogers (1983) states that the jump measure $\Pi$ has a
completely monotone density if and only if the Wiener-Hopf factors
of $X$ are mixtures of exponential distributions. So, we have the
following lemma.
\begin{lemma}
If the jump measure $\Pi$ has a completely monotone density, then
the functions $\psi_q(x)$ and $\psi(x)$ are completely monotone in
$(0,\infty)$.
\end{lemma}

The second result is on Cram\'er's estimate for ruin probability.
Roughly speaking, under suitable conditions, the probability of ruin
decays exponentially when the initial capital becomes larger. The
following lemma is extracted from Kyprianou (2006, Theorem 7.6); see
also Bertion and Doney (1994, Theorem).

\begin{lemma} Assume that $X$ is a L\'evy process which does not
have monotone paths, for which

\noindent (i) (\ref{con-eq1}) holds;

\noindent (ii) there exists a $R>0$ such that $\Psi(-R)=0$ where
$\Psi$ is given by (\ref{intro-eq1});

\noindent (iii) the support of $\Pi$ is not lattice if
$\Pi(\Bbb{R})<\infty$.

\noindent Then,
$$\lim_{x\to\infty}e^{R x}P(\tau<\infty)=\kappa
(0,0)\left(R\frac{\partial \kappa(0,\beta)}{\partial
\beta}\bigg{\vert}_{\beta=R}\right)^{-1},$$ where the limit is
interpreted to be zero if the derivative on the right-hand side is
infinite. For more details,   see Kyprianou (2006, Theorem 7.6).
\end{lemma}

\setcounter{equation}{0}
\section{ Barrier strategy}\label{barrier}
\setcounter{equation}{0}

In this section, we consider a simple barrier strategy for dividend
payments. Under a barrier strategy, if the controlled surplus
reaches the level $b$, then the overflow will be paid as dividends;
and if the surplus is less than $b$, then no dividends are paid out.
Let $\xi_b=\{L_t^b:t\ge 0\}$ be a barrier strategy and $U^b =
\{U_t^b : t \ge 0\}$ be the corresponding controlled risk process.
Note that $U_t^b = X_t - L_t^b$ and $\xi_b\in \Xi$. Moreover, if
$U^b_0\in[0,b]$, then the strategy $\xi_b$ corresponds to a
reflection of the process $X-b$ at its supremum; and if $t\le \tau_b
=\inf\{t\ge 0: U_t^b\le 0\}$, the process $L^b_t$ can be defined by
$L_0^b=0$ and $L_t^b=\sup_{s\le t}[X_s-b]\vee 0$. Note that $L^b_t$
is increasing, continuous and adapted such that the support of the
measure $dL^b_t$ is contained in the closure of the set
$\{t:U_t^b=b\}$. If $U^b_0=x>b$, $L^b_t$ has a jump at $t=0$ of size
$x-b$ to bring $U^b$ back to the level $b$ and a similar structure
afterward.

Let $V_b(x)$ denote the dividend-value function if a barrier
strategy with level $b$ is applied.  Then,
\begin{equation}
V_b(x)=E\left(\int_0^{\tau^b}e^{-\delta t}dL_t^b|U^b_0=x\right).
\label{barrier-eq1}
\end{equation}
The following result shows that $V_b(x)$ as a function of $x$
satisfies the following integro-differential equations with certain
boundary conditions. Note that Paulsen and Gjessing (1997)
established a similar result for a very general jump-diffusion
process.

\begin{theorem}\label{thrm3-1} Assume that the process $X$ have Laplace exponent (\ref{intro-eq1}) and the infinitesimal generator $\Gamma$ is given by (\ref{intro-eq2}). Let $V_b(x)$ be bounded and twice continuously differentiable on $(0, b)$ with a bounded first derivative and with the understanding that we mean the right-hand derivatives at $x=0$.

\noindent (i) \ If $V_b(x)$ solves
$$\Gamma V_b(x)=\delta V_b(x), \quad  0< x <b,$$
together with the boundary conditions
\begin{eqnarray*}
V_b(x) &=& 0, \quad x < 0,\\
V_b(0) &=& 0,  \quad \text{\rm if} \ \ \sigma^2>0, \\
V_b'(b) &=& 1, \\
V_b(x) &=& V_b(b) + x-b, \quad x>b,
\end{eqnarray*}
then $V_b(x)$ is given by (\ref{barrier-eq1}).

Furthermore, let $\psi_b(x)$ be bounded and twice continuously
differentiable on $(0,b)$ with a bounded first derivative and with
the understanding that we mean the right-hand derivatives at $x=0$.

\noindent (ii) If $\psi_b(x)$ solves
$$\Gamma\psi_b(x) =0, \quad  0 < x< b,$$
together with the boundary conditions
\begin{eqnarray*}
\psi_b(x) &=& 1, \quad  x < 0, \\
\psi_b(0) &=& 1, \quad \text{\rm if} \ \sigma^2>0, \\
\psi'_b(b) &=& 0,
\end{eqnarray*}
then $\psi_b(x)=P_x(\tau_b<\infty)$.
\end{theorem}

\noindent {\bf Proof.} \ The  proof of (i) can be done by using the
arguments used in Paulsen and Gjessing (1997). If $\sigma^2>0$, the
process starting from 0 immediately has a negative value.  Hence,
$V_b(0)=0$. Applying Ito's formula to $e^{-\delta (t\wedge
\tau^b)}V_b(U^b_{t\wedge \tau^b})$ gives
\begin{eqnarray*}
e^{-\delta (t\wedge \tau^b)}V_b(U^b_{t\wedge \tau^b})&=&V_b(x)+\int_0^{t\wedge \tau^b} e^{-\delta s}(\Gamma-\delta)V_b (U^b_s)ds\\
&& \quad -\int_0^{t\wedge \tau^b} e^{-\delta s}V_b'
(U^b_s)dL^b_s+M_t,
\end{eqnarray*}
where $M_t$ is a martingale.  Since $V_b(U^b_{\tau^b})=0$ and the
support of the measure $dL^b_t$ is contained in the closure of the
set $\{t:U_t^b=b\}$, taking expectation on both sides of the
equality above yields
\begin{equation}
E\left(e^{-\delta (t\wedge \tau^b)}V_b(U^b_{t\wedge
\tau^b})\right)=V_b(x)-E\left(\int_0^{t\wedge \tau^b} e^{-\delta
s}dL^b_s\right).\label{barrier-eq2}
\end{equation}
Therefore, (i) follows by letting $t\to \infty$ in
(\ref{barrier-eq2}).  The proof of (ii) is entirely analogous to the
proof of (i). \hfill $\Box$

The integro-differential equation
$$\Gamma h(x)=\delta h(x), \quad x>0,$$
has, apart from a constant factor, a unique nonnegative solution
$h(x)$. This together with Theorem 3.1 gives
\begin{equation}
  V_b(x)=\left\{
  \begin{array}{ll}
    \frac{h(x)}{h'(b)},&0\le x\le b,\\
    x-a+\frac{h(b)}{h'(b)}, &x>b.
  \end{array}
  \right.\label{barrier-eq3}
\end{equation}
In particular, if $\Pi\{(0,\infty)\}=0$, then $h$ becomes the
$\delta$-scale function $W^{(\delta)}$ and (\ref{barrier-eq3})
reduces to Proposition 1 of Avram et al. (2007), which was proved by
the excursion theory. See also Renaud and Zhou (2007) and Zhou
(2005) for an alternative approach.

\setcounter{equation}{0}
\section{ Main results and proofs}\label{main}
Define a barrier level by
$$b^*=\sup\{b\ge 0: h'(b)\le h'(x) \ \ \text{for all} \ \ x\ge 0\},$$
where $h'(0)$ is understood to be the right-hand derivative at $0$,
and is not necessarily finite.

We now present the main results of the paper which concerns the
optimal barrier strategy $ \xi_{b^*}$ for a general L\'evy
processes. This  is a continuation of the work of Avram et al.
(2007),  Loeffen (2008), and Kyprianou et al. (2010) in which only
the spectrally negative L\'evy process was considered.

\begin{theorem}\label{thrm4-1} Suppose that $h$ belongs to $C^1 (0,\infty)$ if
 $\sigma=0$ and $\int_{-\infty}^{\infty}\min\{|x|,1\}\Pi(dx)<\infty$ and otherwise belongs to $C^2(0,\infty)$.
If  $\lim_{x\to\infty}h'(x)=\infty$ and $h(x)$ is convex on
$[b^*,\infty)$. Then, the barrier strategy at $b^*$ is an optimal
strategy.
\end{theorem}

For simplicity, we write the L\'evy density $\pi$ as
$$\pi(x)=\left\{\begin{array}{ll} \pi_1(x), &  x>0,\\
 \pi_2(-x), & x<0, \end{array}\right.$$
where $\pi_1, \pi_2$ are L\'evy measures concentrated on
$(0,\infty)$.

\begin{theorem}\label{thrm4-2}
If $\pi_1$ and $\pi_2$ are completely monotone on $(0,\infty)$, then
the barrier strategy at $b^*$ is an optimal strategy.
\end{theorem}

Before proving the main results, we present several lemmas which are
similar to those for spectrally negative L\'evy process. For
$\delta>0$, we consider the following second order
integro-differential equation:
\begin{equation}
\frac{1} {2}\sigma^2 h''(x)+a h'(x)
+\int_{-\infty}^{\infty}[h(x+y)-h(x)-h'(x)y\text{\bf
1}_{\{|y|<1\}}]\Pi(dy)=\delta h(x),\; x>0. \label{main-eq1}
\end{equation}
Set
$$\rho(\delta)=\sup\{\theta: \Psi(\theta)=\delta\}.$$
We assume that $\rho(\delta)>0$. For such a $\rho(\delta)$, we
denote by $P^{\rho(\delta)}$ the exponential tilting  of the measure
$P$ with Radom-Nikodym derivative
$$\frac{dP^{\rho(\delta)}}{dP}\bigg{\vert}_{\Bbb{F}_t}=e^{\rho(\delta)
X(t)-\delta t}.$$ Under the measure $P^{\rho(\delta)}$, the process
$X$ is still a  L\'evy process with Laplace exponent
$\psi_{\rho(\delta)}$ given by
$$\psi_{\rho(\delta)}(\eta)=\psi(\eta+\rho(\delta))-\delta.$$
Let the process $\tilde{X}$ has the L\'evy triplet $(\tilde{a},
\tilde{\sigma}^2, \tilde{\Pi})$ where $\tilde{\sigma}^2=\sigma^2$,
$\tilde{\Pi}(dx)=\tilde{\pi}(x)dx=e^{\rho(\delta) x}\pi(x)dx,$ and
$$\tilde{a}=a+\sigma^2\rho(\delta)+\int_{-\infty}^{\infty}
(e^{\rho(\delta) y}-1)y\text{\bf 1}_{\{|y|\le 1\}}\pi(y)dy.$$
Moreover, \begin{equation} \int_{|x|\ge 1}e^{\rho(\delta)
x}\pi(x)dx<\infty.  \label{main-eq2}
\end{equation}
We refer the reader to Kyprianou (2006) for related discussions.

Note that the law of $\tilde{X}$ is $X$ under the new probability
measure.  Let $\tilde{\psi}(x)$ be the ruin probability for the
L\'evy process $\tilde{X}$. Then, we have the following result.

\begin{lemma}\label{lemma4-1}
The solutions of equation   (\ref{main-eq1}) are proportional to the
function $(1-\tilde{\psi}(x))e^{\rho(\delta) x}$.
\end{lemma}

\noindent{\bf Proof.} \ We first claim that
$E(\tilde{X}_1)=\tilde{\Psi}'(0)$ is always positive where
$\tilde{\Psi}$ is the L\'evy exponent of $\tilde{X}$. In fact,
\begin{eqnarray*}
\tilde{\Psi}'(0)&=&\tilde{a} +\int_{-\infty}^{-1}x\tilde{\pi}dx  +\int^{\infty}_{1}x\tilde{\pi}dx\\
&=&a+\rho(\delta)\sigma^2+\int_{-\infty}^{\infty} ye^{\rho(\delta) y}\pi(y)dy-\int^1_{-1}y\pi(y)dy\\
&=&{\Psi}'(\rho(\delta))>0,
\end{eqnarray*}
since $\Psi$ is strictly convex and increasing on
$[0,\rho(\delta)]$. If $\pi$ is completely monotone, then it follows
from Definition 2.1 that $\pi$ has the form (\ref{con-eq2}). By
using the integrability condition (\ref{main-eq2}), we get
$\int_0^{\rho(\delta)}\mu_{+}(dt)=0$. Therefore,
 $$
\tilde{\pi}(x)\equiv e^{\rho(\delta) x}\pi(x)=\text{\bf
1}_{\{x>0\}}\int_{(\rho(\delta),\infty)}e^{-(t-\rho(\delta))x}\mu_{+}(dt)
+\text{\bf
1}_{\{x<0\}}\int_{(0,\infty)}e^{(t+\rho(\delta))x}\mu_{-}(dt),
$$
which shows that $\tilde{\pi}$ is  completely monotone. So, we see
from Lemma 2.1 that $\tilde{\psi}$ is also completely monotone, and
hence $\tilde{\psi}\in C^{\infty}(0,\infty)$. Furthermore, it
follows from Theorem 3.1 (ii) that $\tilde{\psi}$ solves
$$\tilde {\Gamma}\tilde{\psi}(x) =0,\;\; x>0,$$
together with the boundary conditions
\[
\tilde{\psi}(x) = 1, \quad x < 0,  \qquad \text{and} \qquad
\tilde{\psi}(0) = 1 \quad \text{\rm if} \quad \sigma^2>0,
\]
where $\tilde {\Gamma}$ is the infinitesimal generator for
$\tilde{X}$ given by
$$\tilde{\Gamma}g(x)=\frac{1} {2}\tilde{\sigma}^2 g''(x)+\tilde{a} g'(x)
+\int_{-\infty}^{\infty}[g(x+y)-g(x)-g'(x)y\text{\bf
1}_{\{|y|<1\}}]\tilde{\pi}(y)dy.$$ Then, it can be shown by
straightforward calculations that the function
$$h(x)=(1- \tilde{\psi}(x))e^{\rho(\delta) x}$$
satisfies integro-differential equation (\ref{main-eq1}). Hence, the
result follows. \hfill $\Box$

Moreover, using arguments similar to those in  Avram et al. (2007)
and Loeffen (2008), we have the following two results.

 \begin{lemma}\label{lemma4-2} Suppose that  $h$ belongs to $C^1 (0,\infty)$ if  $\sigma=0$ and $\int_{-\infty}^{\infty}\min\{|x|,1\}\Pi(dx)<\infty$, and otherwise belongs to $C^2(0,\infty)$. If $h$ is  also convex on $[b^*,\infty)$, then
$$(\Gamma-\alpha)V_{b^*}(x)\le 0, \qquad \text{\rm for} \quad x>b^*.$$
\end{lemma}

\begin{lemma}\label{lemma4-3} (Verification lemma) Suppose that $\xi$ is an admissible dividend strategy such that $V_{ \xi}$ is twice continuously differentiable and for all $x>0$
$$\max\{(\Gamma -\alpha)V_{ \xi}(x), 1-V_{ \xi}'(x)\}\le 0.$$
Then, $V_{\xi}(x)=V_{*}(x)$ for all $x$.
\end{lemma}

To end the section, we give the proofs of Theorems 4.1 and 4.2.

\noindent {\bf Proof of Theorem 4.1.} \ The condition
$\lim_{x\to\infty}h'(x)=\infty$ implies that $b^*<\infty$. Clearly,
it follows from the definition of $V_{b^*}$ and  Lemma 4.2 that for
$x>0$
\begin{eqnarray*}
(\Gamma-\alpha) V_{{b^*}}(x) &\le& 0,\\
(1-V_{{b^*}}'(x))(\Gamma-\alpha) V_{{b^*}}(x)&=&0,\\
V_{{b^*}}'(x) &\ge& 1,
\end{eqnarray*}
which in turn imply that for all $x>0$
$$\max\{(\Gamma -\alpha) V_{{b^*}}(x), 1- V_{{b^*}}'(x)\}=0.$$
Hence, the result is a direct consequence of Lemma 4.3. \hfill
$\Box$

\noindent {\bf Proof of Theorem 4.2.} \ If $\pi$ is completely
monotone, then $\tilde{\psi}$ is also completely monotone. Hence,
$\tilde{\psi}(x)$ admits the following representation
\begin{equation}
\tilde{\psi}(x)=\int_0^{\infty}e^{-s x}\mu(ds), \label{main-eq3}
\end{equation} where
$\mu$ is a Borel  measure on $[0,\infty)$. Making use of Cram\'er's
estimate for ruin probability (see Lemma 2.2), we have
$$
\tilde{\psi}(x)\sim  C e^{-(R+\rho(\delta)) x}, \; x\to\infty,
$$
where  $R$ is a positive constant such that $\Psi(-R)=0$, and $C$ is
a nonnegative constant.  Consequently, we have
\[
\lim_{x\to\infty}e^{\rho(\delta) x}\tilde{\psi}(x)=0.
\]
This together with (\ref{main-eq3}) give
$\int_0^{\rho(\delta)}\mu(ds)=0$, and hence
$$e^{\rho(\delta) x}\tilde{\psi}(x)=\int_{\rho(\delta)}^{\infty}e^{-s x+\rho(\delta)
x}\mu(ds)=\int_0^{\infty}e^{-tx}\mu(\rho(\delta)+dt),$$ which is
completely monotone. In particular, we have $(e^{\rho(\delta)
x}\tilde{\psi}(x))'''\le 0$. Thus,
$$h'''(x)=\rho(\delta)^3e^{\rho(\delta) x}- (e^{\rho(\delta) x}\tilde{\psi}(x))'''>0, \quad x>0.$$
That is, $h'$ is strictly convex on $(0,\infty)$. The result follows
from Theorem 4.1 since $\lim_{x\to\infty}h'(x)=\infty$. \hfill
$\Box$

\setcounter{equation}{0}
\section{Examples}\label{exam}

In this section, we present an example with mixed-exponential
jump-diffusion process and a list of completely monotone L\'evy
densities.

\noindent {\bf Example 5.1.} (Mixed-exponential jump-diffusion
process - see also Asmussen et al. (2004) and Mordecki (2004))

Consider the process $X=\{X_t : t\ge 0\}$ given by
\begin{equation}
X_t=at+\sigma B_t+\sum_{k=1}^{N_t}Y_k,  \label{exam-eq1}
\end{equation}
where $B=\{B_t : t\ge 0\}$ is a standard Brownian motion, $N=\{N_t :
t\ge 0\}$ is a Poisson process with parameter $\lambda$, $Y=\{Y_k :
k\ge 1\}$ is a sequence of independent and identically distributed
random variables with density
\begin{equation}
\pi(x)=\left\{\begin{array}{ll} p \sum_{j=1}^m A_j\eta_j e^{-\eta_j x}, &  x>0,\\
q\sum_{j=1}^n B_j\varsigma_j e^{\varsigma_j x}, &
x<0,\end{array}\right.  \label{exam-eq2}
\end{equation}
with $p, q \ge 0$, $p+q=1$, $ \eta_j, \varsigma_j>0$, $ A_j, B_j\ge
0$, $\sum_{j=1}^m A_j=1, and \sum_{j=1}^n B_j=1.$ As usual, we
assume that the processes $B$, $N$, and $Y$ are independent. When
$m=n=1$, the process $X$ reduces to the double-exponential
jump-diffusion process (see Kou and Wang (2003)). Since $\pi$
defined in
  (\ref{exam-eq2}) is completely monotone, a direct application of Theorem 4.2
shows that the barrier strategy at $b^*$ is an optimal strategy.
\hfill $\Box$

Besides Example 5.1, other examples of L\'evy processes with
completely monotone densities can be found in the literature.
Several of them are listed below (for details, see Bagnoli and
Bergstrom (2005), Loeffen (2008), Yin and Wang (2009) and Kyprianou,
Rivero and Song (2010)):

\begin{itemize}
\item $\alpha$-stable process with L\'evy density:
$\pi(x)=\lambda x^{-1-\alpha},\; x>0$ with $\lambda>0$ and
$\alpha\in (0,1)\cup (1,2);$

\item one-sided tempered stable process (particular cases include gamma process ($\alpha=0$) and inverse Gaussian process
($\alpha=1/2$)) with L\'evy density: $\pi(x)=\lambda
x^{-1-\alpha}e^{-\beta x},\; x>0$ with $\beta, \lambda>0$ and $-1\le
\alpha <2;$

\item the associated parent process with L\'evy density:
$\pi(x)=\lambda_1 x^{-1-\alpha}e^{-\beta x}+\lambda_2
x^{-2-\alpha}e^{-\beta x},$ $ x>0$ with $\lambda_1, \lambda_2>0$ and
$-1\le \alpha <1.$
\end{itemize}
Note that they all satisfy the condition
$\int_0^{\infty}\pi(x)dx=\infty$. Moreover, some distributions with
completely monotone density functions are given below:

\begin{itemize}
\item Weibull distribution with density: $f(x)=cr x^{r-1}\text{e}^{-c x^r}, \;x>0,$ with $c>0$ and $0<r<1$;

\item Pareto distribution with density: $f(x)=\alpha (1+x)^{-\alpha-1}, \;x>0,$ with $\alpha>0$;

\item mixture of exponential densities: $f(x)=\sum_{i=1}^n A_i \beta_i \text{e}^{-\beta_i x}, \;x>0,$ with $A_i>0, \beta_i>0$ for
$i=1,2\cdots,n$, and $\sum_{i=1}^n A_i=1$;

\item gamma distribution with density:
$$f(x)=\frac{x^{c-1}e^{-x/\beta}}{\Gamma (c)\beta^c}, \quad x>0,$$
with $\beta>0$ and $0<c\le 1$.
\end{itemize}

\noindent{\bf Remark.} \ For a spectrally negative  L\'evy process,
whenever the L\'evy measure has a density which is log-convex, then
Kyprianou et al. (2010) showed that the optimal strategy is a
barrier strategy. Since there is no condition on the upward jumps in
Theorem 4.2, it might be possible to generalize the result of
Kyprianou et al. (2010) to the present situation. That is, for any
L\'evy process with arbitrary positive jumps, if the L\'evy density
of negative jumps is log-convex, then the optimal dividend strategy
is a barrier strategy. However, we are not able to give a formal
proof of such a conjecture. \hfill $\Box$

\vskip0.3cm

\noindent{\bf Acknowledgements}

\; We would like to thank the two anonymous referees who gave us
many constructive suggestions and valuable comments on the previous
version of this paper. The research of Kam C. Yuen was supported by
a university research grant of the University of Hong Kong. The
research of Chuancun Yin was supported by the National Natural
Science Foundation of China (No.10771119) and the Research Fund for
the Doctoral Program of Higher Education of China (No.
20093705110002).

\end{document}